\documentclass[aps,rmp]{revtex4}
\pdfoutput=1

\usepackage{amsfonts,amssymb,amsmath,bm}
\usepackage{graphics,epsfig}
\usepackage{verbatim,ulem}
\usepackage{hyperref}
\usepackage{breakurl}
\usepackage{here}
\usepackage{color}
\usepackage{comment}
\usepackage[ngerman,french,english]{babel}

\def\be{\begin{eqnarray}}
\def\en{\end{eqnarray}}
\def\UB{\begin{equation}}
\def\UE#1{\label{#1}\end{equation}}
\def\BE{\begin{equation}}
\def\EE#1{\label{#1}\end{equation}}

\def\nL{\nabla^{\rm L}}
\def\rf#1{(\ref{#1})}

\def\ba{{\bm as}}
\def\bu{{\bm u}}
\def\bx{{\bm x}}
\def\bX{{\bm X}}

\def\ba{{\bm a}}
\def\bF{{\bm F}}
\def\bn{{\bm n}}
\def\bV{{\bm V}}
\def\bom{{\bm\omega}}

\hypersetup{
     colorlinks   = true,
     citecolor    = blue,
     linkcolor=blue
}

\bgroup\catcode`\#=12\egroup \newcommand{\HD}[2]{\href{#1}{{\color{blue} #2}}}

\begin{document}
\title{A contemporary look at Hermann Hankel's 1861 pioneering\\ ~~~~~~~~~~~~~~~~~work on
  Lagrangian fluid dynamics}
\author{Uriel Frisch}
\affiliation{Universit\'e C\^ote d'Azur,~OCA,~Lab.~Lagrange,~CS~34229,~06304~Nice~Cedex~4,~France}
\author{G\'erard Grimberg}
\affiliation{Instituto de Matematica - UFRJ, Av. Athos da Silveira Ramos 149, CT - Bloco
C, Cid. Univ. - Ilha do Fund\~ ao. Caixa Postal 68530 21941-909 Rio de Janeiro
- RJ - Brasil}
\author{Barbara Villone}
\affiliation{INAF, Osservatorio Astrofisico di Torino, Via Osservatorio, 20, 10025 Pino Torinese, Italy}
\date{\today}
\begin{abstract}
   The present paper is a companion to the paper by Villone and Rampf
   (2017), titled ``Hermann Hankel's
\textit{On the general theory of motion of fluids},
an essay including an English translation of the complete Preisschrift
from 1861'' together with connected documents. Here we
   give a critical assessment of Hankel's work, which covers many
   important aspects of fluid dynamics considered from a
   Lagrangian-coordinates point of view: variational formulation in
  the spirit of Hamilton for elastic (barotropic) fluids, transport
   (we would now say Lie transport) of vorticity, the Lagrangian
   significance of Clebsch variables, etc.  Hankel's work is also put 
   in the perspective of previous and future work. Hence, the action
   spans about two centuries: from Lagrange's 1760--1761 Turin paper on
   variational approaches to mechanics and fluid mechanics problems to
   Arnold's 1966 founding paper on the geometrical/variational  formulation of
   incompressible flow. The 22-year old Hankel --- who was to die 12
   years later ---  emerges as a highly
   innovative master of mathematical fluid dynamics, fully deserving
   Riemann's assessment that his  \textit{Preisschrift} contains ``all
   manner of good things.''

\end{abstract}
\maketitle

\section{Introduction}
\label{s:intro}

It has been known for over two centuries that there are two ways for describing the motion of a fluid.
In the \textit{Eulerian approach} the velocity, the density and the pressure
 are expressed in terms of the coordinates $\bx =(x,y,z)$ of the point of
measurement and of the time $t$, whereas in the \textit{Lagrangian approach} they are
expressed in terms of the coordinates $\ba =(a,b,c)$ of the initial position
of the fluid particle and of the time elapsed. In modernized notation, the two systems are related by the
Lagrangian map $\ba \mapsto \bx(\ba,t)$, which satisfies the characteristic
equation $\dot\bx = \bu (\bx,t)$, where $\bu$ is the velocity and the dot
denotes the Lagrangian (or material) time
derivative. Before the Eulerian/Lagrangian terminology was in
use, in his famous 1788 treatise ``M\'echanique analitique'', Lagrange
warned his readers that ``The  [Eulerian] equations are much simpler
than
the [Lagrangian] ones; thus, in the theory of fluids,  one should
prefer
the former over the latter.'' Even a couple of centuries after ``Lagrange's
curse,'' Lagrangian coordinates have been occasionally
viewed as an ``agony''.\footnote{For Lagrangian coordinates, see
  Euler, \hyperlink{Euler1760}{1760}; Lagrange,
  \hyperlink{Lagrange1788}{1788}. For Lagrange's curse, see Lagrange,
  \hyperlink{Lagrange1788}{1788}:453 and also p.~450.  For the ``agony'', see 
  Price, \hyperlink{Price}{2006}:31.}

Nevertheless, in the 20th century, Lagrangian coordinates have been used extensively
by mathematicians, in analyzing the well-posedness of ideal
incompressible fluid
flow, governed by the Euler equations.The goal was to prove theorems
regarding the existence, regularity and uniqueness of the solution to
the 
initial-value problem, where one prescribes the flow at time $t=0$. In the
20s and the 30s, the works of G\"unther, Lichtenstein and Wolibner
were the first to show that ideal incompressible fluid flow with suitable initial
regularity, stays so for at least a finite time (in three dimensions)
and for all times in two dimensions (in a bounded domain).  These
results have barely been improved since and the situation is not much
better for viscous flow, at least in three
dimensions.\footnote{G\"unther, \hyperlink{Guenther}{1926}. Lichtenstein,
  \hyperlink{Lichtenstein}{1927}. Wolibner,
  \hyperlink{Wolibner}{1933}.}

A major reason why Lagrangian approaches work better than Eulerian
ones for such questions has been identified several decades
later by Ebin and
Marsden, using an (infinite-dimensional) geometrical approach to the
dynamics of ideal incompressible fluids, pioneered by Arnold. The key
is that the application of Lagrangian time derivatives do not lose
spatial derivatives, whereas Eulerian time derivatives do lose spatial
derivatives.
What this means,  is that if the initial
velocity field has spatial derivatives up to order $p\ge 2$, then 
initial Eulerian time derivatives exist only up to the same order
$p$, whereas Lagrangian time derivatives exist up to any
order. Roughly, the reason is that  in an Eulerian formulation
a time derivative $\partial_t$ is always accompanied by a space derivative
$\bu\cdot\nabla$, whereas in a Lagrangian formulation these two terms
combine to form the Lagrangian time derivative.\footnote{Ebin and
  Marsden, \hyperlink{Ebin}{1970}. Arnold, \hyperlink{Arnold1966}{1966}; Arnold \& Khesin,
\hyperlink{ArnoldKhesin1998}{1998}.}

We mention that there has been for several decades a strong interest
for Lagrangian approaches in geophysical fluid
dynamics, cosmology and magnetohydrodynamics.\footnote{Frisch and  Villone, \hyperlink{Frisch}{2014} and references therein.}

In the 19th century there was relatively little use of Lagrangian
coordinates, not so much because of Lagrange's curse but because there
was a more objective reason  to stay
away from Lagrangian coordinates: After the work on the Navier--Stokes equations, it was realized that for
most realistic flows one cannot ignore the role of viscosity and to
achieve this, it is indeed easier to use Eulerian coordinates.

Some key theoretical results obtained during the 19th century were
actually formulated as Lagrangian conservation theorems, but their
derivations were given in
Eulerian coordinates. This includes   Helmholtz's result
on the conservation of vorticity flux and Thomson/Kelvin's result on
the conservation of circulation.\footnote{Helmholtz, \hyperlink{Helmholtz1858}{1858}; Thomson (Lord Kelvin), \hyperlink{Thomson}{1869}}

But there are also 19th century studies that make end-to-end use of
Lagrangian variables. Two stand out: Cauchy's derivation of the three-dimensional invariants, now
called the Cauchy invariants and the results by Hankel presented here
 and in the companion paper. Both Cauchy's and Hankel's were prized contributions
on a prescribed theme. For Cauchy it was ``the problem of waves on the
surface of a liquid of arbitrary depth'', so that his results on
Lagrangian invariants were somewhat incidental and occupied only a small
fraction of his long memoir. For Hankel the situation was quite different:
he was being asked to write a memoir on the theory of fluids centered on the use of Lagrangian coordinates,
and that is basically what he did.\footnote{For Cauchy and how his work on
  invariants was mostly forgotten in the 19th century, see Cauchy,
  \hyperlink{Cauchy1815and1827}{1815/1827}; Frisch  and Villone,
  \hyperlink{Frisch}{2014} \S\S\, 2--4.}

Around the end of the 1850s, among the mathematicians and physicists at the University
of G\"ottingen, such as Weber, Dirichlet
and  Riemann, all of
whom were interested in fluid mechanics, awareness of the
importance of Lagrangian coordinates started rising. Following  the death
of Dirichlet in 1859, an extraordinary prize was
set up, which was won in 1861 by Hermann Hankel,  a  student in
Mathematics on leave from
Leipzig University. We shall come back to the detailed circumstances around that
prize in Section~\ref{s:context}, but we wish to explain briefly why
two companion papers in this issue of EPJ~H are dedicated to
Hankel's Preisschrift (Prize manuscript), the present one and another
one containing the full translation of the Preisschrift with  related documents.\footnote{For the Preisschrift, see Hankel, \hyperlink{Hankel1861}{1861}. For the companion paper, see Villone and Rampf, \hyperlink{Villone}{2017}.} 

A detailed study of Hankel's Preisschrift reveals indeed a 
deep understanding of Lagrangian fluid mechanics, with variational
methods and differential geometry present throughout. This did not
escape Riemann's attention who, in the official assessment for the special
prize committee set up by the faculty, wrote that the manuscript
contained 'all manner of good things'. The style of writing of Hankel may
occasionally appear awkward, but is actually fully consistent with
mid-19th century German scientific writing tradition.\footnote{For 
  exchanges between Riemann and Weber relating to Hankel's application
  for the prize and for the official assessment, see the companion paper by Villone \& Rampf, \hyperlink{Villone}{2017}.}

Later Hankel did make
important contributions to mathematics until he died at age 34.  His
work on  hydrodynamics did not achieve much visibility. In the 19th
century,
Hankel's hydrodynamics work is cited several times by
Beltrami, and Auerbach covers the work rather extensively. 
In addition, there are references concerned exclusively with Hankel's
proof of the Stokes theorem, namely Lipschitz and
Lamb.\footnote{Beltrami,  \hyperlink{Beltrami1871}{1871}. Auerbach,
  \hyperlink{Auerbach1881}{1881}. Lipschitz,
  \hyperlink{Lipschitz}{1870-1871}. Lamb, \hyperlink{Lamb1895}{1895}.}

  In the 20th century, we have to wait for Truesdell to find extensive
  coverage of Hankel's hydrodynamics, which is also mentioned by
  Darrigol.\footnote{Truesdell, \hyperlink{Truesdell1954a}{1954a};
    Truesdell and Toupin, \hyperlink{TruesdellToupin1960}{1960}.
    Darrigol, \hyperlink{Darrigol}{2005}.}

The paper is organized as follows. Section~\ref{s:context} describes the
context of the mathematical and fluid dynamical research in Germany and
particularly in G\"ottingen, at the time of Hankel's Preisschrift.

Section~\ref{s:state} describes the tools and methods available at the
time and also their relation to more modern methods used currently, e.g., in
differential geometry.

Section~\ref{s:variational} is about Hankel's first variational approach 
to elastic (barotropic) fluid dynamics in the spirit of Hamilton's least action principle.

Section~\ref{s:cauchy} shows how Hankel derived the Helmholtz theorem for the
conservation of vorticity flux by using Cauchy's 1815 Lagrangian reformulation
of the incompressible Euler equations. The conservation of circulation plays
there a role as an intermediate step.\footnote{Cauchy, \hyperlink{Cauchy1815and1827}{1815/1827}.}

Section~\ref{s:clebsch} is about Hankel's rederivation of the Clebsch-variable
approach to the Euler equations using clever switches between Lagrangian and
Eulerian coordinates, now called pullbacks and pushforwards. This  is
followed  by two subsections: Section~\ref{ss:1859} briefly describes how
Clebsch introduced the variables known by his name two years before Hankel's
work; Section~\ref{ss:topology} is about a topological pitfall of Clebsch 
variables, understood only in the last fifty years.

Concluding remarks are presented in Section~\ref{s:conclusions}, in particular
about the difficult birth and rebirth of variational methods.

\section{Context of Hankel's 1861 prize}
\label{s:context}

Three important scientists connected with the Georg--August university
of G\"ottingen 
around the middle of the 19th century played a major role in the
establishment of the prize that was attributed to Hankel. They are
the physicist Wilhelm Eduard Weber (1804--1891), particularly known
for his work on electricity and magnetism and closely associated to one of the
brightest stars of G\"ottingen, Carl Friedrich Gauss (1777--1855), and
two mathematicians, Gustav Peter Lejeune-Dirichlet (1805--1859) and 
Bernhard Riemann (1826--1866), who successively held the chair of
Gauss (Riemann had also been assistant to Weber).\footnote{Biographies
  of Weber, Dirichlet and Riemann may be found in the Dictionary of
  Scientific Biography, \hyperlink{dictionary}{1970-1980}. Biographies
  of Weber, Dirichlet and Riemann: these are by  Woodruff, \hyperlink{Woodruff}{1976}, Ore, \hyperlink{Ore}{1859} and  Freudenthal, \hyperlink{Freudenthal}{1975}, respectively.}

Around 1850 the Georg--August university was quickly turning into a major
center for mathematics and physics. Gauss was of course still an important
driving element but until the political liberalisation that followed the 1848
upheavals in many European countries, Georg--August suffered from serious lack
of freedom for academics, which in turn led to a dramatic drop in student
enrollment. For example, Weber was one of the seven professors (including the
Brothers Grimm) who were summarily dismissed in 1837, although he was not
banned from G\"ottingen and was eventually reinstated in 1849, becoming one of
its most influential professors.  Weber played a leading role in establishing
the prize discussed here. Dirichlet had been teaching for many years in
Berlin, where he and Carl Gustav Jacob Jacobi (1804--1851) attracted many
bright students, among whom Riemann. After the death of Gauss in 1855, with
the recommendation of Weber, Dirichlet was appointed to Gauss's chair in
G\"ottingen. He died however in 1859 and the chair went to Riemann, who in
December 1859 presented to the Academy of G\"ottingen an unpublished work of
Dirichlet on hydrodynamics and another unpublished work of his own on
hydrodynamics.\footnote{For the G\"ottingen seven and its impact on student
    enrollment, see Hunger, \hyperlink{Hunger}{2002}. 
    For Weber's role in the hiring of Dirichlet,
    see Minkowski, \hyperlink{Minkowski}{1905}. 
    For the presentation of Dirichlet's hydrodynamics unpublished work, 
    see Nachrichten, \hyperlink{NachLD}{1857}, \hyperlink{NachR}{1859} and for
    the posthumous published work, see Dirichlet, \hyperlink{LejeuneDirichlet}{1860}.  
    For Riemann's 1859
    presentation of his own 
    hydrodynamics work, see Nachrichten, \hyperlink{NachR}{1859} 
    and the published work, Riemann, \hyperlink{Riemann}{1861b}.}

The death of Dirichlet and his posthumous paper  ``On a problem of
hydrodynamics,'' edited by Richard Dedekind (1831--1916)\footnote{Lejeune-Dirichlet, \hyperlink{LejeuneDirichlet}{1860}.}, which
discusses the advantages and disadvantaged of using Eulerian or
Lagrangian coordinates, are
highlighted in the official announcement as one of the  motivations of the prize for which Hankel
applied. For over one century, establishing and studying the equations of  hydrodynamics had been considered a
major challenge by mathematicians.  In 1858 Hermann Helmholtz
(1821--1894) had discovered the basic laws of vortex dynamics and
developed an interesting analogy with the nascent electromagnetic
theory: velocity is expressible in terms of vorticity using the same
Biot--Savart law that relates magnetic fields to currents. Helmholtz's
theorem on the conservation of vorticity flux through a surface
element following the flow, did probably strike Riemann as a  potentially
Lagrangian result, although Helmholtz gave an Eulerian
derivation. Riemann had himself been interested in the Eulerian
vs. Lagrangian approaches to hydrodynamics and observed that they both
originated with Euler. It is thus not surprizing that the extraordinary
prize which was established after these events
requested a fully Lagrangian derivation of Helmholtz's
result.\footnote{Riemann's views on the origin of Eulerian and
  Lagrangian coordinates are discussed in \S\,1  of the Preisschrift. Helmholtz, \hyperlink{Helmholtz1858}{1858}.}

Here is the official statement of the Georg--August university about
the prize \textit{The most useful equations for determining fluid motion may be presented
in two ways, one of  which  is Eulerian, the other one is
Lagrangian.
The illustrious  Dirichlet pointed out in the posthumous unpublished
paper
``On a problem of hydrodynamics'' the almost completely overlooked
advantages  of the
Lagrangian way, but he was  prevented from unfolding this way further  by a
fatal illness. So, this institution asks for a theory of fluid motion  based on
the equations of Lagrange, yielding, at least, the laws of vortex motion
already derived in another way by the illustrious
Helmholtz}.\footnote{The original statement was in Latin: ``Aequationes generales
motui fluidorum determinando inservientes
duobus modis exhiberi possunt, quorum alter Eulero, alter Lagrangio
debetur. Lagrangiani modi utilitates adhuc fere penitus neglecti
clarissimus Dirichlet indicavit in commentatione postuma 'de problemate
quodam hydrodynamico' inscripta; sed ab explicatione earum uberiore
morbo supremo impeditus esse videtur. Itaque postulat ordo theoriam
motus fluidorum aequationibus Lagrangianis superstructam eamque eo
saltem perductam, ut leges motus rotatorii a clarissimo  Helmholtz
alio modo erutae inde redundent.''}
The prize was set up on 4th June 1860 and had a
deadline of end of March 1861.

Enters Hermann Hankel, who arrived in the spring 1860 in G\"ottingen,
as a 21-year old mathematics student, coming   from the University of  Leipzig. In G\"ottingen, Hankel worked on the prize, which he won
as single contestant. This was  a common situation for such thematic
prizes, but the prize was not necessarily awarded. According to Zahn, while in G\"ottingen, Hankel attended the  lectures of Riemann  during
the 1860 summer term and the 1860--1861 winter term on such topics as: ``The
mathematical theory of  gravity, electricity  and magnetism'' and
``The theory of partial differential equations with  applications  to
physical problems.''   Hankel stayed in G\"ottingen until
the Fall of 1861, defending his doctoral dissertation
in Leipzig during the same year. Then he left for Berlin to extend
his mathematical training, attending 
lectures by   Weierstrass and others.\footnote{According to Zahn, \hyperlink{Zahn}{1874}:584, Hankel came to G\"ottingen in April
  1860 to study with Riemann and others. For
  the catalogue of lectures in G\"ottingen, 
    see Nachrichten, \hyperlink{NachH}{1860a}, \hyperlink{NachLancement}{1860b}.}

\section{State of the art near the middle of the 19th century, as used
by Hankel}
\label{s:state}

The most popular way of formulating the equations of fluid dynamics
and, more generally of continuous media, from the time of Lagrange to
the middle of the 19th century and beyond was the method of
\textit{virtual work} (also known as ``virtual displacements'' or
``virtual velocities''), applied to the various external forces and
the inertial force. A considerable extension and synthesis of previous work by
Johann (John)
Bernoulli (Principle of virtual velocities) and of D'Alembert's Principle,
Lagrange presented it as a founding principle of
mechanics.\footnote{For the virtual velocities Principle, see
  Bernouilli, Johann (John), \hyperlink{Bernoulli1717}{1717}; for
  D'Alembert's Principle and its discussion, see  D'Alembert,
  \hyperlink{Dalembert}{1743};  Lagrange,
  \hyperlink{Lagrange1788}{1788} and Truesdell,
  \hyperlink{Truesdell}{1960}:188--190.} In modernized notation this
amounts roughly to the following. One avoids Euler's, traditional way
of writing Newton's equation for each infinitesimal element
(``molecule''), namely $\bF = m\ba$, where $\bF$, $m$ and $\ba$ are
the force, the mass and the acceleration, respectively. (Henceforth we
use modernized notation, mostly boldface vector notation.) Instead,
two changes are performed. First, one lumps together into a total external
force $\bF_{\rm ext}(\bx,t)$ all the external forces acting on the molecule
near $\bx$ at time $t$, including the ``inertial force''
$-m\ba$. Second, one gives a weak formulation (with space derivatives
taken in the sense of distributions) of $\bF_{\rm ext}(\bx)=0$, by
equating to zero, at any given time, the total work $\int \bF_{\rm
  ext}(\bx,t)\cdot \delta \bx(\bx,t) d^3\bx$ of these forces under assumed
infinitesimal virtual displacements $\delta{\bx}(\bx,t)$.  An instance
is (1) of \S\,4 of HT.\footnote{Henceforth, HT designates the
  translation of Hankel's Preisschrift in the companion paper by
  Villone and Rampf, \hyperlink{Villone}{2017}. By (n) \S\,p we
  understand the equation in \S\,p that was given the number (n) by
  Hankel (such numbers are preserved in the translation). Since many
  of Hankel's equations are unnumbered, in the translation all the
  equations have been numbered, using the format [p.n], which means
  the nth equation of \S\,p.} As was common throughout the second half
of the 18th century and the 19th, Hankel systematically considers both
the case of incompressible fluids and that of elastic fluids having a
well-defined functional relation between the pressure and the
density. Nowadays, the latter would be called barotropic.  In the
incompressibe case, the infinitesimal virtual displacements are chosen
to preserve the volume, i.e., $\nabla \cdot \delta\bx =0$. An instance
is [3.2] in HT.
  
When reading Hankel, it must be kept in mind that he always starts from the
virtual work formulation  and then derives what we generally
call Euler's equations. So, actually he goes from a weak to a strong
formulation. Since no analytic convergence issues were considered
at the time of Hankel, this does not matter.

Interestingly, Hankel did not shy off the now standard notation of
Euler, avoiding virtual displacements. However it was probably
important for the young man, still a student, to demonstrate that he
was fully familiar with the notation most common at that time. 

A second important point when reading Hankel's Preisschrift, is the central
role played by the duality of approaches mentioned in the theme of the prize,
namely the use of either Eulerian coordinates. referred to by Hankel as
``first manner/method'' or Lagrangian coordinates, referred to by Hankel as
``second manner/method''. 

A central feature of Hankel's Preisschrift is the constant juggling with
Eulerian and Lagrangian coordinates. Hankel performs changes of
variables from Eulerian to Lagrangian coordinates (or vice-versa) on
both scalar functions and on differential forms (in modern geometrical
parlance these are called pullbacks and  pushforwards). As long as he works with scalar functions, these are just
changes of coordinates. As was common at the time, he does not change the name
of the function so transformed; neither will we. When he performs  pushbacks and
pullforwards on differential forms (or on gradients), the changes in
the 
differentials are controlled by  the Jacobian matrix $\partial x_i/\partial a_j$. Because of the
heavy use of virtual work (weak) formalism, Hankel frequently needs the 
volume element. Denoting the mass density  by $\rho(\bx,t)$  and the initial
density by $\rho_0 (\ba)$, one has of course
\UB \rho(\bx,t)\, d^3 \bx = \rho_0(\ba)\,d^3\ba.\UE{masscons}
Although he makes frequent indirect use of this mass conservation relation, the closest 
he comes to writing it, is in HT [2.9] which, in modernized form, reads
\UB \det \left(\frac{\partial x_i}{\partial a_j}\right)  = \frac{\rho_0}{\rho}.\UE{massconsweak}

Hankel stresses in
 \S\,1 that, according to  his
advisor Bernhard Riemann, Euler is the originator of both
approaches, the Eulerian and the Lagrangian ones. This issue has been
revisited by Truesdell.  Of course, Hankel had to make
central use of Lagrangian coordinates in the Preisschrift. Actually, he succeeded
amazingly well in this, and particularly so because he became fully aware
of Cauchy's 1815 work on Lagrangian coordinates.\footnote{Truesdell
  \hyperlink{Truesdell1954b}{1954b}. Cauchy, \hyperlink{Cauchy1815and1827}{1815/1827}.}

Finally, Hankel and most theoretical mechanicians of his time were
quite familiar with developments of variational methods that started
about one century earlier with the work of Euler and Lagrange. It is not
our purpose here to review the history of variational methods.
Lagrange, a few years after he started discussing with Euler his ideas on
analytical methods for variational problems, tried applying such methods
to various mechanical problems, including fluid dynamics. This attempt was
not very successful, in particular because Lagrange used the Maupertuis form of
the action, based just on the kinetic energy; furthermore, he used
``Lagrangian'' coordinates only to monitor the evolution of volume elements.\footnote{For the history of variational methods, 
see Goldstine, \hyperlink{Goldstine}{1980}; Fraser, \hyperlink{Fraser2003}{2003}. For the early variational work, see Lagrange \hyperlink{Lagrange1760}{1760--1761}; for
``Lagrangian'' coordinates, see his \S\,XLIV, p.~276.}

A crucial change of perspective occurred in the 1830th when Hamilton
realized that for most mechanical systems one should include in the action 
not only a kinetic but also a potential contribution and that it can be
advantageous to switch to so-called canonical variables. This line of 
research was pursued and disseminated widely all over Europe by Jacobi
and many others.\footnote{Hamilton, \hyperlink{Hamilton1834}{1834}; \hyperlink{Hamilton1835}{1835}. Jacobi, \hyperlink{Jacobi1837}{1837}; \hyperlink{Jacobi1866}{1866}.} 

However such development concerned chiefly solid mechanics. Fluid mechanical
developments had to wait until around 1860 with the work of Helmholtz and
Clebsch and, of course Hankel's Preisschrift, as we shall now
see.\footnote{Helmholtz, \hyperlink{Helmholtz1858}{1858}. Clebsch, \hyperlink{Clebsch1857}{1857}; \hyperlink{Clebsch1859}{1859}.}

\section{Hankel's first variational formulation for elastic fluids}
\label{s:variational}
Hankel's derivation of a variational formulation for the dynamics
of elastic/barotropic flow is presented at the beginning of his \S\,5,
where it takes barely more than one page. The corresponding Eulerian equations (HT
\S\,4 (4) [4.5]) are written in his previous section. In modernized form
they read (after pulling out a  factor $\rho$ in the  first equation):
\UB\rho\left[\partial_t \bu +\bu\cdot\nabla\bu -\bX+\frac{1}{\rho}\nabla p\right] = 0,\qquad
\partial_t\rho+\nabla\cdot(\rho\bu)=0,\UE{eulerelastic}
where $\bu$ is the velocity, $p$ the pressure, $\rho$ the density,
$\nabla$ the Eulerian spatial gradient and
$\bX$ the external force per unit volume, assumed by Hankel to derive from a potential ($\bX = \nabla V$), as is for example the case for gravitation.
This equation has an equivalent virtual-work (weak) form also derived in Hankel's
previous section (HT \S\,4 (1) [4.1])
\UB\int d^3\bx\left\{\rho\left[(\bX -\ddot\bx)\cdot \delta\bx\right]
+p\,\nabla\cdot\delta\bx\right\}=0,\UE{weakeulerelastic}
where $\ddot\bx$ denotes the acceleration of the fluid particle, that
is the Lagrangian/material time derivative of the velocity $\dot\bx$.

From \eqref{weakeulerelastic}, using $\bX = \nabla V$ and integrating
by parts the term involving the pressure, Hankel obtains, after a global
change of sign, [5.2]
\UB\int d^3\bx\,\rho\left[\ddot\bx -\nabla V +\frac{1}{\rho} \nabla p\right]\cdot \delta\bx=0.\UE{elasticIbyP}
Then, using the barotropic relation $\rho = \phi(p)$ and introducing the
function $f(p) = \int dp/\phi(p)$, Hankel observes that $\frac{1}{\rho}
\nabla p = \nabla f(p)$. Next come two important steps, not particularly
underlined by Hankel. First, introducing the new function 
$\Omega \equiv V-f(p)$ he makes use of $\delta \Omega=\nabla \Omega\cdot \delta\bx$. Second,
he performs a pullback to Lagrangian coordinates, using \eqref{masscons}, to
obtain [5.6]:
\UB\int d^3\ba\,\rho_0(\ba)\left[\ddot\bx \cdot \delta\bx
  -\delta\Omega\right]=0.\UE{5.6}
It is now obvious --- both for us  and for a mechanician of Hankel's time ---
that the variational equation \eqref{5.6}, which should hold at any time
in, say, the interval $[0,T]$, is equivalent, after a time integration by parts, to the vanishing of the variation
of an action,  namely (HT \S\,5 ((1) [5.9], essentially)
\UB \delta\,\int d^3\ba\,\rho_0(\ba)\int_0^T dt\, \left[\frac{1}{2} |\dot
\bx|^2+ \Omega\right]=0,\UE{hankelaction}
with the usual constraint that the infinitesimal variations $\delta \bx$ vanish
at time $0$ and $T$. Here we have modernized Hankel's
$\left(\frac{ds}{dt}\right)^2$ (where $s$ is the arclength) into
$|\dot\bx|^2$. Actually, the form with $ds^2$ is quite useful, as Hankel
will explain a little further.

Hankel has thus given a variational formulation  for the elastic/barotropic
fluid. It involves an action that is the space-time integral of the
Lagrangian density function $\rho_0(\ba)\left[\frac{1}{2} |\dot
\bx|^2+ \Omega\right]$, namely the difference of the kinetic energy $(1/2)\rho_0(\ba)|\dot \bx|^2$ and of the potential energy $-\rho_0(\ba)\Omega$, the latter being itself the sum of  the potential
energy $-\rho_0(\ba) V$ due to external forces and of the elastic potential 
$\rho_0(\ba) f(p)$.

Here Hankel pauses to comment on the nature of his result. He qualifies
it as a theorem possessing some  analogy to the principle of least action. 
Of course, the latter was historically based on the Maupertuis action, involving
just the kinetic energy. Only with Hamilton and Jacobi did mechanicians
realize that the simplest  correct variational formulation for conservative
systems must use an action involving both the kinetic and the
potential energy. Such modern variants of the least action principle were
commonplace for discrete mechanical system but not yet for fluids. One 
possible exception is the work of Clebsch on incompressible flow, to which
we shall return in Section~\ref{s:clebsch}.\footnote{For the Maupertuis
 action, see Maupertuis, \hyperlink{Maupertuis1740}{1740}, \hyperlink{Maupertuis1744}{1744}; Euler, \hyperlink{Euler1744}{1744}; \hyperlink{Euler1750}{1750}. Hamilton, \hyperlink{Hamilton1834}{1834}; \hyperlink{Hamilton1835}{1835}. Jacobi, \hyperlink{Jacobi1866}{1866}.
Clebsch, \hyperlink{Clebsch1857}{1857}; \hyperlink{Clebsch1859}{1859}.}  

At this point, Hankel makes an important observation about the use of 
his variational formulation. One essential ingredient in the action appearing
in \eqref{hankelaction} is the square of the line element $ds ^2$. In
Cartesian coordinates we have $ds ^2 =dx^2 +dy^2+dz^2$. In another
coordinate system, involving arbitrary curvilinear coordinates, denoted $\rho_1$,
$\rho_2$ and $\rho_3$, the squared line element will be some other quadratic 
form of the $d\rho_i$ ($i=1,2,3$) and the volume element can be easily
evaluated in arbitrary coordinates. Hankel observes that it is much
easier to perform the change of variables on the action integral than doing
it directly on the equation of motion. He actually devotes 13
out of the 53 pages of the
Preisschrift to deriving such equations in a number of different coordinate
systems. He then points out that there is no need to stay with three variables and
that the same process can be carried in an $n$-dimensional space. Here, he
was probably  influenced by Riemann's teaching, who in the same year
1861 submitted (unsuccessfully) to the French Academy a prize essay on
the conduction of heat  using the same notation for curvilinear
coordinates. Nevertheless, Hankel does of course
 not try leaving ordinary Euclidean space in favor of hydrodynamics in  curved Riemann
spaces.\footnote{Riemann, \hyperlink{Riemann1861a}{1861a}.}

\section{From Cauchy to Helmholtz (via the circulation theorem)}
\label{s:cauchy}
The \S\S\,6--9 of Hankel present a fully Lagrangian-coordinates-based
proof of one of Helmholtz's main results, that the flux of the  
vorticity
vector $\nabla \times \bu$ through an infinitesimal surface element
remains
unchanged as we follow the surface element in its Lagrangian
motion. Helmholtz proved this for a three-dimensional incompressible
flow
subject to a potential external force. Hankel's proof is not just
Lagrangian,
as required by the rules of the prize, instead of  
incompressibility, he assumes that the flow is
elastic/barotropic. Yet, as he points out,
the essence of his proof applies equally well to an incompressible
flow.
The starting point in \S\,6 is the pullback to Lagrangian coordinates
of the elastic/barotropic Euler momentum equation \eqref{eulerelastic} (HT [6.1])
\UB\ddot x_k \frac{\partial x_k}{\partial a_i}
    -\frac{\partial\Omega}{\partial a_i}=0,\UE{eulerpullback}
where $\Omega \equiv V -f(p)$ as before and repeated indices are
summed upon. This is obtained as follows: the
momentum equation \eqref{eulerelastic} may be rewritten either as $\ddot \bx - \nabla
\Omega =0$ or, in the language of differential geometry,  as  $\ddot
  \bx\cdot d\bx -d\Omega=0$. The  Eulerian gradient $\nabla \Omega$ is then converted into a Lagrangian
  gradient $\nL \Omega$ through multiplication by the Jacobian matrix
  $\nL \bx$  of the Lagrangian map. In the language of differential geometry,
this is the pullback of the vanishing $1$-form $\ddot
  \bx\cdot d\bx -d\Omega$ from
  Eulerian to Lagrangian coordinates. Henceforth we shall rewrite
  equations such as \eqref{eulerpullback} in more compact notation as
  \UB\ddot x_k \nL x_k = \nL\Omega.\UE{eulerpullbacksv}.

  So far, Hankel has been essentially following Lagrange's own
  rewriting of the hydrodynamical equations in ``Lagrangian
  coordinates''. The next step will be to follow Cauchy in his 1815
  Prize work. Hankel takes a Lagrangian curl of
  \eqref{eulerpullbacksv} to kill the gradient on the r.h.s. and
  obtain (HT [6.2]) \UB \nL \times \left(\ddot x_k \nL x_k\right) =0. \UE{precauchy}
  Still, following exactly Cauchy, Hankel  observes that this equation can be integrated in time to yield
  \UB \nL \dot x_k \times \nL x_k = \bom_0,\UE{cauchy} where $\bom_0
  \equiv \nL \times \bu_0$ is the initial vorticity, which obviously
  does not depend on time. Hence, the three components of the
  l.h.s. are time-independent. They are nowadays called the Cauchy
  invariants.\footnote{For the pullback to Lagrangian coordinates, see
    Lagrange, \hyperlink{Lagrange1788}{1788}:446.  For the Cauchy invariants, see Cauchy, \hyperlink{Cauchy1815and1827}{1815/1827};
    Frisch and  Villone, \hyperlink{Frisch}{2014}. It is noteworthy
    that Cauchy's derivation of equations \eqref{precauchy}-\eqref{cauchy}
    has been apparently  rediscovered (in words only)  by Dedekind in
    his own footnote 1 to Dirichlet's posthumous paper
    (Lejeune-Dirichlet, \hyperlink{LejeuneDirichlet}{1860}), where Dedekind
    also states  that Dirichlet cannot have failed to notice this, as
    well as the consequence, nowadays called the Cauchy vorticity formula.}

The next step is Hankel's own: he rewrites \eqref{cauchy} as \UB \nL
\times\left(\dot x_k \nL x_k\right) = \bom_0,\UE{hankelcauchy} thus pulling
out a curl operator. This is crucial because it will allow him to use what is
now called the Stokes theorem, which as is known, relates the line integral along a
closed contour of a  vector field to the flux of its curl accross a surface
bounded by the contour. The proof  is presented in Hankel's
\S\,7 and constitutes the first published proof of the theorem. This
theorem was actually known much before Hankel's 1861 work, as explained
by Katz. Furthermore, Hankel benefited from close association with
Riemann who had proved
a two-dimensional version of the theorem, itself  related to an
earlier resulf of Cauchy. It is interesting that Riemann's work
discussed the topogical aspects of the theorem, in particular issues of
simply- and multiply-connected  domains, whereas Hankel assumed not only a 
simply-connected
domain but also that the surface can be resolved uniquely for
the $z$-coordinate in terms of the $x-$ and $y$-coordinates. These
simplifications somewhat reduce the scope of Hankel's proved theorem,
but did not matter, since eventually, Hankel only needed to work with
an infinitesimal surface, as assumed by Helmholtz.\footnote{Katz,
  \hyperlink{Katz1979}{1979}. Cauchy,
  \hyperlink{Cauchy1846}{1846}. Riemann, \hyperlink{Riemann}{1857};
  Roy, \hyperlink{Roy2011}{2011}:362.}

The derivation of the Helmholtz theorem, in \S\S\,8--9, proceeds as follows. 
In the Lagrangian space, Hankel considers  a connected
(\textit{zusammenh\"angende}) surface, here denoted by $S_0$, whose
boundary is a curve, here denoted by $C_0$, and applies the Stokes
theorem to \eqref{hankelcauchy}.
 Hankel then infers that the flux through $S_0$ of the 
r.h.s. of \rf{hankelcauchy}, namely the initial vorticity,
is given by the circulation along $C_0$ of the initial velocity $\bu_0$ (HT
[8.5] and also
\S\,8 (2) [8.6]):
\UB \int_{C_0}\bu_0\cdot d\ba=\int_{S_0}\bom_0
\cdot \bn_0 \,d\sigma_0,\UE{hankelhtu}
where $\bn_0$ denotes the local unit  normal to $S_0$ and $d\sigma_0$ the surface
element.
 Then, he similarly handles the l.h.s. of \rf{hankelcauchy}
and first notices that  (HT [9.3])
\UB u_k\nL x_k \cdot d\ba= \bu\cdot d\bx.\UE{hankelndt}
 He thus
obtains the Eulerian circulation, an integral over the curve $C$ where are
presently located  the fluid particles initially on $C_0$ (HT [9.5]):
\UB \int_C \bu\cdot d\bx = \int_{S_0}\bom_0
\cdot \bn_0 \,d\sigma_0.\UE{hankelndu}
Eq.~\rf{hankelndu}, together with \rf{hankelhtu}, although they do not
appear in the same numbered paragraph of Hankel, are clearly a statement of the
standard circulation theorem, mostly associated to the name of
William Thomson (Kelvin). However, as already stated, Hankel's
formulation of the Stokes theorem is restricted to simply-connected domains, which is enough to recover Helmholtz's
flux-invariance result. Hence Hankel's 1861 circulation result is less general
than Kelvin's 1869 result in which multiply-connected domains are taken
into consideration. Of course, Hankel was aware of Riemann's emphasis
on connectedness, since he cited Riemann's 1857 founding paper on the subject.\footnote{Thomson (Lord Kelvin), \hyperlink{Thomson}{1869}. Riemann, \hyperlink{Riemann1857}{1857}.}

At this point, it is clear that Hankel also proved the constancy
in time of the flux of the vorticity through any finite surface moving
with the fluid. Letting  this surface shrink to an
infinitesimal element, he  obtains Helmholtz's theorem.\footnote{Helmholtz,
  \hyperlink{Helmholtz1858}{1858}.} 

Hankel has thus obtained both global and local forms of the Helmholtz
theorem and, incidentally, of the circulation theorem.

\section{Clebsch variables in Hankel's Lagrangian formulation}
\label{s:clebsch}

Alfred Clebsch (1833--1872) entered  K\"onigsberg University in 1850,
attended lectures  by Jacobi's disciples, Otto Hesse, Friedrich
Richelot and Franz Neumann. He graduated in 1854,  moved to Berlin where
he became \textit{Privat Dozent}  in 1858. He went on to  Karlruhe in
1858, to  Gie{\ss}en in 1863 and finally became \textit{Professor} in
G\"ottingen in 1868, where he occupied the chair previously held by
Gauss, Lejeune Dirichlet and Riemann, until his premature death in
1872. His early research work (1857--1859),  focussed on hydrodynamics and the
calculation of variations and was directly inspired by the methodology developed
by Jacobi and Hesse. (Eventually, Clebsch became one of the editors of the
complete works of Jacobi.)  Around 1860 his research interests
shifted towards
the theory of algebraic invariants, where he used the analytical view
of projective geometry implemented by Pl\"ucker (later, Clebsch and
Felix Klein  edited Pl\"ucker's  posthumous works).  After 1863 in
Gie{\ss}en, he worked on Abelian functions with Gordan, following
the path opened by Riemann. After 1868, he became the chief organizer
of
mathematical research in G\"ottingen, creating the journal
\textit{Matematische Annalen} 
and developing contacts with  young French mathematicians, such as
Jordan or Darboux, and with the English school of Sylvester and Cayley. He
interacted with Felix Klein and Sophus Lie. Here, we are only
interested in the very early period of Clebsch's life before Hankel's 1861
Preisschrift. More detailed biographies of Clebsch are of course
available.\footnote{For detailed biographical information on Clebsch,
see Burau, \hyperlink{Burau}{1970--1980} and Various Authors, \hyperlink{VarAuth1873a}{1873a}; \hyperlink{VarAuth1873b}{1873b}.}

The Clebsch-variable formulation of three-dimensional  incompressible flow is
nowadays usually written as follows. There are two three-dimensional
time-dependent scalar fields $\phi(\bx,t)$ and $\psi(\bx,t)$ that are
material invariants, that is, they satisfy
\UB (\partial_t+\bu\cdot \nabla)\phi=0 \quad \text{and}\quad 
(\partial_t+\bu\cdot \nabla)\psi=0.\UE{phipsimaterial}
Furthermore, there is a scalar field $F(\bx,t)$, such that  the velocity field $\bu$ is coupled back to the fields $\phi$ and
$\psi$ by 
\UB \bu = \nabla F + \phi\,\nabla \psi\quad \text{and}\quad \nabla\cdot \bu =0.
\UE{pseudopfaff}
Actually, Clebsch had a somewhat complicated proof that a  velocity field
thus constructed satisfies the incompressible Euler equations. We shall come
back  to this a little further, in Section~\ref{ss:1859}.\footnote{Clebsch, \hyperlink{Clebsch1859}{1859}.} 

In \S\,10--11 of
the Preisschrift, Hankel gives a very short proof of this result, while
stressing that the equations just given are  written in Eulerian
coordinates. Hankel's
proof uses his modified form of the Cauchy invariants equation
\eqref{hankelcauchy} and performs a reversion from Lagrangian to Eulerian
coordinates, which one now calls a pushforward.

First, he assumes that the \textit{initial} vorticity $\bom_0(\ba)\equiv
\bom(\ba,0)$ has the following Pfaff-type representation, obtained by
taking at $t=0$ the curl of \eqref{pseudopfaff} (HT \S\,10 (2) [10.4]): \UB \bom_0 = \nL \times\bu_0 = \nL
\phi_0 \times \nL \psi_0 = \nL \times \left(\phi_0\,\nL
\psi_o\right),\UE{darbouxinitial} where we recall that $\nL$ denotes the
Lagrangian gradient.   We will come back to the validity of this 
representation  at the end of Section~\ref{ss:1859}.\footnote{Pfaff, \hyperlink{Pfaff}{1814}.}

The functions here denoted by $\phi_0$ and $\psi_0$ were written
by Hankel without the subscript zero but he stressed that this assumption is
made only at time $t=0$.
Hankel then substitutes \eqref{darbouxinitial} in the r.h.s. of the Cauchy
invariants equation \eqref{hankelcauchy}. He notices that there are now
Lagrangian curls on both sides, which he can remove provided he adds  the 
Lagrangian gradient of a suitable function $F$. He thus obtains:
\UB \dot x_k\nL x_k = \nL F +\phi_0\,\nL \psi_0. \UE{haha}

Now comes a crucial step. Hankel reverts to Eulerian coordinates and
also switches from Lagrangian gradients (of $F$ and $\psi_0$) to Eulerian
ones.
For the latter step he must use what is now called the inverse Jacobian
matrix of the Lagrangian map, that is $\nabla \ba$. On the l.h.s. there
is already a direct Jacobian matrix $\nL \bx$. In modern language these
two matrices, inverses of each other,  multiply to the identity, so that after this 
multiplication  the l.h.s of \eqref{haha} becomes just $\dot \bx = \bu$ 
and \eqref{haha} becomes exactly \eqref{pseudopfaff}, hence concluding
Hankel's derivation.

These manipulations, described here in matrix language,  are essentially what 
Hankel does when invoking his equations \S\,6 (2) (HT [6.4])  and \S\,2 (3) 
(HT [2.14]). Furthermore, since he performs a change of variables, he might
have changed the name of his functions. Instead, following the 
usage of the time, he resorts to warning language: 
``in so far as time $t$ appears explicitly or implicitly through $x$, $y$,
$z$.'' We avoid this problem by switching from $\phi_0$ and $\psi_0$ to
$\phi$ and $\psi$, although they are actually the same functions in Lagrangian
coordinates, being material invariants.

It is of interest to point out that, shortly after the beginning of this
clever and simple proof, Hankel makes a digression concerning his geometrical
Helmholtzian vision of Clebsch variables.
Knowing, from Helmholtz, that vortex
lines are material, he observes that under the assumption that $\bom_0 = \nL
\phi_0 \times \nL \psi_0$ a simple way to ensure the materiality of vortex
lines is to take the surfaces $\phi = \text{Const.}$ and $\psi =
\text{Const.}$ to be material, namely to have $\dot \phi =0$ and $\dot \psi =
0$.  Hankel also adds, at the end of \S\,11, after completion of the proof,
that Clebsch himself did not bring out the meaning of his variables. Probably
Hankel meant the geometrical meaning. Connected to this, there is also an
important footnote [A.26] 
to which we shall come back after having
summarized the Clebsch approach.

\subsection{The Clebsch 1859 approach}
\label{ss:1859}

Our purpose here is to find out how the 1859 work of Clebsch may have
influenced the subsequent work of Hankel, who was fully aware of it and
cites it several times. For the convenience of the reader we shall use 
Hankel's notation (except for a certain Hamiltonian function present
only in Clebsch). Since Clebsch was including the pressure $p$ in the 
potential of external forces, we shall use $V-p$ where Clebsch uses $\bV$. 
For simplicity we assume that the constant and uniform density is unity.
Furthermore we shall specialize to the three-dimensional incompressible
case, commenting afterwards on the generalization by Clebsch.
Here, the presentation does not aim at analyzing Clebsch's
hydrodynamical work  exhaustively. It would be of interest to do so and
also to translate his 1857 and 1859 papers.\footnote{Clebsch, \hyperlink{Clebsch1857}{1857}; \hyperlink{Clebsch1859}{1859}.}

The starting point are the Euler equations
\UB \nabla (V-p) = \partial_t \bu +\bu\cdot \nabla \bu;\quad \nabla\cdot \bu =0,\UE{eulerclebsch}
the first of which is written as a 1-form differential relation, with the infinitesimal
element denoted $\delta \bx$
\UB \delta(V-p) = \delta \bx \cdot\left(\partial_t \bu +\bu\cdot \nabla \bu\right).\UE{deltaeuler}
Clebsch assumes that, \textit{at any time} and not just initially, the
velocity has the (Clebsch) representation \eqref{pseudopfaff}. After various
algebraic
manipulations, he arrives at the following relation
\UB \delta\left\{ V -p -(1/2)|\bu|^2 -\partial_t F  -\phi\,\partial_t
\psi\right\} = \dot \phi \,\delta \psi -\dot \psi \,\delta \phi, \UE{clebsch8}
where $\dot \phi \equiv (\partial_t +\bu\cdot\nabla) \phi$ and $\dot \psi
\equiv (\partial_t +\bu\cdot\nabla) \psi$. 

Clebsch then notes that the l.h.s., being the variation of a scalar function,
so must the r.h.s. By counting the number of independent variations he argues
that the r.h.s.  is of the form $\delta \,\Pi$ where $\Pi$ can only depend on $\phi$,
$\psi$ and the time $t$. From there, immediately follows that we have a pair
of canonical Hamilton equations
\UB \dot \phi =-\frac{\partial\Pi}{\partial \psi}, \quad \dot \psi
=\frac{\partial\Pi}{\partial \phi}.\UE{hamiltontodie}
However, Clebsch shows that by performing a suitable canonical transformation,
this Hamiltonian $\Pi$ morphs into a vanishing one. The details of this
transformation, which actually originated with Jacobi are now
standard.\footnote{Jacobi, \hyperlink{Jacobi1890}{1836-1837/1890}: 393; Lanczos, \hyperlink{Lanczos1970}{1970}: Chap.~8, Sec.~2, p.~238; Landau and Lifshitz,  \hyperlink{Landau1976}{1976}:148.} As a consequence, there is a pair of new Clebsch
variables, which we shall here still denote by $\phi$ and $\psi$  whose
material time derivatives vanish, so that they satisfy 
\eqref{phipsimaterial}.

Clebsch does not comment on the vortical and Lagrangian interpretation
of his variables, but his work appeared only shortly after the
Helmholtz 1858 vorticity paper, which is not cited. Hankel's
subsequent work is of course much simpler, mostly because he makes
use of Lagrangian coordinates, as required for the
Prize.\footnote{Helmholtz, \hyperlink{Helmholtz1858}{1858}.}

Here it is appropriate to mention Hankel's own assessment of some of the 1859
work of Clebsch, as it appears in the first part of the  footnote [A.26] at
the end of \S\,11.\footnote{Footnote  1 on p.~45 of the German
  original.}  Hankel stresses that, on the one hand, the Clebsch derivation
(probably of \eqref {clebsch8}) can be further simplified and, on the other
hand, that his direct proof --- that the variables $\phi$ and $\psi$ can be
chosen to be material --- may be adapted to a purely Eulerian
proof. Eventually,  such a purely Eulerian derivation of the material
character of the Clebsch variables will be given.\footnote{For
  Eulerian derivations, see Duhem,
  \hyperlink{Duhem1901}{1901} and references therein.}

Of course there are also modern derivations of such results, using the
tools of differential geometry and in particular the concept of
Lie-advection invariance. For example, proceeding as in Besse
and Frisch,  one can take $\phi$ and
$\chi$ to be two Lie-advection-invariant $0$-forms (modern parlance for
stating that they are material-invariant scalars). Then, the exterior
derivative ${\mathrm d} \psi$ and the product $\phi{\mathrm d} \psi$
will be Lie-advection-invariant $1$-forms. The latter's exterior
derivative ${\mathrm d} \phi\wedge {\mathrm d} \psi$ is then a
Lie-advection-invariant $2$-form. If initially this $2$-form is taken
to be the vorticity $2$-form, it will stay so.  Thus the vorticity
vector (the Hodge dual of the vorticity $2$-form) can be written at all time as
$\nabla \phi \times \nabla \chi$. Except for the modernized language,
this is essentially the proof given by Hankel in \S\S\,10-11.  Hankel here makes use
of concepts and calculation procedures, which will take their full
significance only within the framework of Lie's theoy. 
\footnote{Besse and Frisch, \hyperlink{BesseFrisch2017}{2017}.}

A variational formulation of the incompressible Euler equations, in the spirit
of Hamilton and Jacobi, can then be given as follows. One observes that, $\phi$
and $\psi$ being now material invariants, \eqref{clebsch8}, after integration
over space and time (over an interval $[0,T]$ at the ends of which the
variations are taken to vanish) and after an integration by parts in time,
becomes \UB\delta\int_0^T dt\, d^3\bx\left\{\psi\,\partial_t \phi
-(1/2)|\bu|^2 +V -p\right\} = 0. \UE{clebschhamilton}
Eq.~\eqref{clebschhamilton} constitutes a variational formulation for the
canonically conjugate variables $\psi$ and $\phi$ and a Hamiltonian $H =
(1/2)|\bu|^2 -V +p$, provided $V$ and $p$ are somehow reexpressed functionally
in terms of $\psi$, $\phi$ and time.

We should also mention that Clebsch did not restrict his work to three
dimensions. He was working in dimension $2n+1$ where $n$ is an arbitrary
positive integer and he had not one but $n$ pairs of variables.

\subsection{Clebsch variables and global topology}
\label{ss:topology}

There is  a  topological issue that could not be perceived at the time
of Clebsch, Hankel and all their contemporaries. 
Indeed,  about a century later, with
the work on Kolmogorov-Arnold-Moser (KAM) theory and the work relating
helicity to vortex line knotedness, it became clear that writing the vorticity in
the Pfaffian way, namely as a vector product of two gradients, as in  \eqref{darbouxinitial}, 
may not hold globally in space. Hadamard already had pointed out that
the Clebsch representation is generally only local.
 
In the three-dimensional case, it was shown by Darboux that a Pfaff
representation can be made only locally, near a point of non-vanishing
vorticity. To understand why such a representation can generically not be
extended globally, we note that the vortex lines are at the intersections of
the surfaces of constant $\phi$ and constant $\psi$, as observed by
Hankel. Such intersections constitute typically closed lines. However the most
general topology of the integral lines of a three-dimensional vector field of
vanishing divergence is given by the KAM theorem and has chaotic topology. One
well-known instance are the ABC flows.  A global Clebsch representation holds
only in special cases, one instance being the Taylor--Green flow.  This
difficulty with the usual Clebsch variables can however be bypassed by using
multi-Clebsch variables, that is more than one pair of Clebsch variables, even
in three dimensions, where Clebsch had just one pair.\footnote{For the Pfaff
  problem, see Pfaff, \hyperlink{Pfaff}{1814}. Clebsch,
  \hyperlink{Clebsch1862}{1862}. Darboux, \hyperlink{Darboux}{1882}.
  Hadamard, \hyperlink{Hadamard1903}{1903}:79--81. For the history of
  the Pfaff problem, see Hawkins, \hyperlink{Hawkins}{2013}:155--204. 
  For the KAM theory, see
  Kolmogov, \hyperlink{Kolmogorov}{1954}; Arnold, \hyperlink{Arnold1963}{1963}; 
  Moser, \hyperlink{Moser}{1962}.
  For the ABC flow, see Dombre,  Frisch, Greene, and H\'enon, \hyperlink{Dombre}{1986} and references
  therein. For the Taylor--Green flow, see Taylor  and Green, \hyperlink{Taylor}{1937};
  Brachet,  Meiron, Orszag,  Nickel, Morf and Frisch, \hyperlink{Brachet1983}{1983}. 
  For multi-Clebsch variables, see Zakharov, L'vov  and
  Falkovich, \hyperlink{Zkaharov}{1992}:28.}

There is however another mechanism that can prevent a global Clebsch
representation: if the vortex line are closed but knotted. Such flows are said
to possess helicity.\footnote{Moreau, \hyperlink{Moreau}{1961};  Moffatt, \hyperlink{Moffatt}{1969}; Bretherton, \hyperlink{Bretherton}{1970};
  Kotiuga, \hyperlink{Kotiuga}{1991}.}

\section{Conclusions}
\label{s:conclusions}

Before attempting to give a modern assessment of Hankel's
Preisschrift, it is of interest to look at Hankel's own assessment of
his work, as it appears in his 4-page summary published in
Fortschritte der Physik (Progress in Physics).\footnote{Hankel, \hyperlink{Hankel1863}{1863}. This is
signed with a somewhat cryptic ``Hl'', which is however explained in the
list of authors at the end of the volume.} He first presents his variational formulation and observes that it is a very
convenient way for studying hydrodynamics in an arbitrary system of coordinates
by just changing the expression of the element of length $ds^2$. This shows
his high regard for Riemann. He then turns to Cauchy's Lagrangian formulation and applications to the
Helmholtz vortical motion. He also states Riemann's view that the Cauchy
invariants are an expression of conservation of rotation. Hankel also points
out that the exploitation of Cauchy's result requires the use of what we now
call the Stokes theorem, which Hankel had just rediscovered independently and
of which it is not clear that Helmholtz was aware in 1858. 
Hankel then gives the essence of the derivation of
the Helmholtz theorem from the Cauchy invariants. Finally, Hankel gives his
own view on Clebsch variables, in which  their material invariance is established
directly by using Lagrangian coordinates.
More than 150 years later, we do not feel that this summary has aged
much. We would just add Hankel's establishement of the conservation of
circulation, albeit he never emphasizes this result, now
considered by many as one of the most important of ideal fluid dynamics.

We have left for these concluding remarks the thorny issue of how the
variational formulation of the incompressible Euler equations has
appeared. Given the important role it has played in the work of the
last fifty years (initiated by Arnold's 1966 SDiff formulation), we
have felt that this is a worthwhile question.\footnote{Arnold, \hyperlink{Arnold1966}{1966}.}

Already in his founding paper on three-dimensional incompressible
flow, Euler speculated that the equation may have a least action
formulation in the sense of Maupertuis. The Maupertuisian action is a
sum or integral of the velocity $v$ times the particle infinitesimal
displacement $ds$. Since $ds = v dt$, this action reduces --- except
for an overall factor two --- to the kinetic energy contribution of
the modern (Hamiltonian) action, leaving out the potential energy
contribution. The peculiarity of an incompressible fluid is that
there is no genuine potential energy: everything is subsumed in the
geometrical constraint of volume conservation.  Hence, a Maupertuisian
attack can be successful, provided incompressibility is correctly
taken into account. Lagrange has precisely introduced the appropriate
tool, the multipliers known by his name.  One would thus expect to
find the first variational/least-action formulation of the
incompressible flow in Lagrange's \textit{opus magnum} of 1788. It is actually
there, but only in an indirect way: Lagrange points out that for a
rather general class of mechanical systems, to which the principle of
live forces is applicable, a least action formulation leads to the
same equations as his virtual velocity formulation.\footnote{Lagrange,
  \hyperlink{Lagrange1788}{1788}:188--189.} Here the principle of live forces is essentially
equivalent to the modern statement of total energy conservation.

Of course, as is well known, around the middle of the 1760s, Lagrange
stopped being very enthusiastic with variational approaches and
switched to virtual velocity/virtual power approaches, which he viewed
as more general. Before this methodological transition,
he published in 1760--1761 a fairly long and rather
difficult paper, where variational methods were applied to a number of
different mechanical problems, including hydrodynamics.  He gave a
Maupertuisian least action derivation of the elastic Euler equations,
using what we would now call the conservation of total energy. Doing a
better job required the development of Hamilton's version of the least
action principle in which both the kinetic and the potential energy
appear in the Hamiltonian and the Lagrangian. In the case of elastic
hydrodynamics, this project was carried out by Hankel in a most elegant
and simple manner using Lagrangian coordinates (see
Section~\ref{s:variational}).\footnote{For Lagrange's 
transition to virtual velocities, see Lagrange, \hyperlink{Lagrange1764}{1764}; Dahan-Dalmenico,
\hyperlink{Dahan}{1992}:39--42. For Lagrange's early work on least action methods applied to
mechanics, see Lagrange, \hyperlink{Lagrange1760}{1760--1761}; Fraser, \hyperlink{Fraser1983}{1983}.}

Why does Hankel not also discuss the case of incompressible flow? Actually he
touches on this issue in the second half of the footnote [A.26] at the end of
\S\,11.  He notes that
in his 1859 paper Clebsch gives an Eulerian variational formulation for
three-dimensional incompressible flow and states that ``In another form, this
is the theorem I presented in \S\,5 [discussed in our
  Section~\ref{s:variational}], deduced immediately from the principle of
virtual velocities.'' Hankel may not have desired or found the time to discuss
an issue already touched by Clebsch. Furthermore, Hankel may have thought, as
we do, that Lagrange essentially was aware of the variational formulation for
incompressible flow, nearly one century before Hankel's work.

It is interesting that it took yet another century
before Arnold's 1966 variational formulation of ideal incompressible
3D flow in terms of SDiff, the space of volume-preserving
diffeomorphisms. 

Arnold's 1966 work  was however much more than just a variational formulation. Lagrange, as he stresses in the
\textit{Avertissement} (Warning) at the beginning 
 of his 1788 Analytic Mechanics, engaged in a much-needed
 degeometrization of mechanics, which brought considerably more rigor
 in the presentation. Eventually --- but this took two centuries ---
 Arnold became the father of the regeometrization of fluid mechanics. The
 central quantity in the geometrical approach of Arnold, of Ebin and
 Marsden and of others, already present in the work of Cauchy and of
 Hankel, is the Lagrangian map. Dealing with such maps requires
 an entirely novel type of geometry,
 capable of handling  infinite-dimensional Riemannian manifolds and their geodesics,
 concepts that where very far from  ripe at the time of Riemann and
 Hankel.\\

\noindent ACKNOWLEDGMENTS 

We are grateful to Nicolas Besse, Olivier Darrigol, Cornelius Rampf, 
Emmanuele Tassi and anonymous referees for useful remarks and suggestions. GG and BV wish to thank the
Observatoire de la C\^ote d'Azur and the Laboratoire J.-L.~Lagrange 
for their hospitality and financial support.

{}
\end{document}